# On the Construction of Doubly Even Order Magic Squares

Grasha Jacob[1], Dr. A. Murugan[2]

[1] Research and Development Centre, Bharathiar University, Coimbatore – 641046, India
[Associate Professor, Dept. of Computer Science, Rani Anna Govt College for Women, Tirunelveli, India]
(e-mail: grasharanjit@gmail.com)
[2] Associate Professor, Dept. of Computer Science, Dr. Ambedkar Govt Arts College, Vyasarpadi, Chennai, India

**Abstract:** *Magic squares have been an enthralling topic in mathematics for centuries. They are formed by filling in all the cells of a square matrix with the numbers starting from one so that the sum of all rows, columns, and diagonals is the same. Magic squares have applications in entertainment, music and even cryptography. This paper focuses on a simple and easy method to construct doubly even magic squares.*

**Keywords:** Square matrix, Magic Square, Doubly even order, Magic Constant

## 1 INTRODUCTION

A magic square of order n is a square matrix or array of $n^2$ numbers such that the sum of the elements of each row and column, as well as the main diagonals, is the same number, called the magic constant denoted by σ(M). Generally, the entries are thought of as the natural numbers 1, 2, ..., $n^2$, where each number is used exactly once. Magic Squares exist for all orders, n ≥ 1 with an exception when n = 2. The case n = 1 is trivial as it consists of a single cell containing the number 1 alone. The constant sum in every row, column and diagonal is called the magic constant or magic sum. The magic constant of a normal magic square depends only on n. For normal magic squares of order n = 3, 4, 5, 6, 7, 8 … the magic constants are 15, 34, 65, 111, 175, 265… respectively. Magic squares are classified into three types: odd, doubly even and singly even. A doubly even magic square is a square matrix of order n, where n is divisible by four. A singly even magic square is a square matrix of order n, where n is even but not divisible by four.

## 2 PROPERTIES OF MAGIC SQUARES

- The magic constant of a normal magic square depends only on n, where n is the order of the magic square. The magic constant is given by the formula, σ(M) = $n(n^2+1)/2$.
  For normal magic squares of order n = 3, 4, 5, 6, 7, 8 … the magic constants are 15, 34, 65, 111, 175, 265… respectively.
- A magic square will remain magic if any number is added to every number of a magic square.
- A magic square will remain magic if any number multiplies every number of a magic square.
- The transpose of a magic square is also a magic square [3].
- A magic square will remain magic if two rows, or columns, equidistant from the center are interchanged.

- An even order magic square will remain magic if the diagonally opposite quadrants are interchanged.
- An even order magic square will remain magic if the columns or the rows are exchanged. If the 1st and the 2nd columns/rows are exchanged and the 3rd and the 4th columns/rows are exchanged, this is called a 1-2-3-4 exchange.

A most-perfect (MP) magic square MM is of doubly-even order (n = 4, 8, . . .) and has the following properties:

- two elements that are n/2 elements apart along all diagonals (including broken ones in both directions) sum to $n^2+1$.
- the elements of all 2 by 2 subsquares (including broken top-bottom, broken left-right, and the four corners) sum to the constant $2n^2+1$.

## 3 MAGIC SQUARES OF DOUBLY EVEN ORDER

A doubly even magic square is constructed by drawing x's through each 4 x 4 subsquare, filling all squares in sequence and then replacing each entry $a_{ij}$ on a crossed-off diagonal by n x n + 1 - $a_{ii}$ or, equivalently, reversing the order of the crossed-out entries. For n = 8, the crossed-out numbers are originally 1, 4, ... , 61, 64, so entry 1 is replaced with 64, 4 with 61, etc as in figure 1.

| 1  | 63 | 62 | 4  | 5  | 59 | 58 | 8  |
|----|----|----|----|----|----|----|----|
| 56 | 10 | 11 | 53 | 52 | 14 | 15 | 49 |
| 48 | 18 | 19 | 45 | 44 | 22 | 23 | 41 |
| 25 | 39 | 38 | 28 | 29 | 35 | 34 | 32 |
| 33 | 31 | 30 | 36 | 37 | 27 | 26 | 40 |
| 24 | 42 | 43 | 21 | 20 | 46 | 47 | 17 |
| 16 | 50 | 51 | 13 | 12 | 54 | 55 | 9  |
| 57 | 7  | 6  | 60 | 61 | 3  | 2  | 64 |

Figure 1. Magic Square of Order 8

## 4 PROPOSED METHOD

There are many ways to construct magic squares, but the standard way is to follow certain formulas for generating the patterns. In this paper, a doubly even magic square is constructed by drawing x's through the main diagonals of each 4 x 4 subsquare, filling all squares in sequence and then replacing each entry $a_{ii}$ other than the crossed-off diagonal by n x n + 1 - $a_{ij}$ and having crossed off entries as such. Figure 2 symbolizes the proposed method.

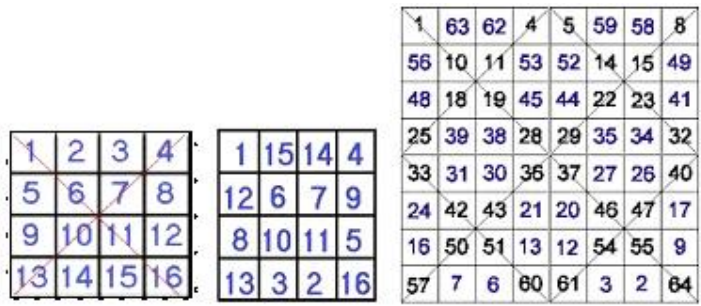

Figure 2. Construction of Magic Squares

## Algorithm Magic Square

Input : Order of matrix, n

Output : Magic Square, A

1. *Read the value of n, where n is doubly even.*
2. *Declare a matrix A of order n X n. Fill all the cells of the matrix with zeros.*
3. *Divide the matrix into blocks of size 4 x 4.*
4. *Fill the diagonal of each 4 X 4 block with the corresponding cell value, $c_{ij}$.*
5. *Trace through each row and fill all the cells of the entire matrix with zero entries with the value of $n \times n - c_{ij}$.*

**End Magic Square**

### 4.1  Steps in the construction of Doubly Even order Magic Square
**Step 1**
Read the value of n, where n is doubly even.
Let us assume n = 8.

**Step 2**
Consider a matrix A of order 8 X 8. Fill all the cells of the matrix with zeros.

| 0 | 0 | 0 | 0 | 0 | 0 | 0 | 0 |
|---|---|---|---|---|---|---|---|
| 0 | 0 | 0 | 0 | 0 | 0 | 0 | 0 |
| 0 | 0 | 0 | 0 | 0 | 0 | 0 | 0 |
| 0 | 0 | 0 | 0 | 0 | 0 | 0 | 0 |
| 0 | 0 | 0 | 0 | 0 | 0 | 0 | 0 |
| 0 | 0 | 0 | 0 | 0 | 0 | 0 | 0 |
| 0 | 0 | 0 | 0 | 0 | 0 | 0 | 0 |
| 0 | 0 | 0 | 0 | 0 | 0 | 0 | 0 |

**Step 3**
Divide the matrix into blocks of size 4 x 4.

**Step 4**
Fill the diagonal of each 4 X 4 block with the corresponding cell value, $c_{ij}$.

| 0 | 0 | 0 | 0 | 0 | 0 | 0 | 0 |
|---|---|---|---|---|---|---|---|
| 0 | 0 | 0 | 0 | 0 | 0 | 0 | 0 |
| 0 | 0 | 0 | 0 | 0 | 0 | 0 | 0 |
| 0 | 0 | 0 | 0 | 0 | 0 | 0 | 0 |
| 0 | 0 | 0 | 0 | 0 | 0 | 0 | 0 |
| 0 | 0 | 0 | 0 | 0 | 0 | 0 | 0 |
| 0 | 0 | 0 | 0 | 0 | 0 | 0 | 0 |
| 0 | 0 | 0 | 0 | 0 | 0 | 0 | 0 |

| 1 | 0 | 0 | 4 | 5 | 0 | 0 | 8 |
|---|---|---|---|---|---|---|---|
| 0 | 10 | 11 | 0 | 0 | 14 | 15 | 0 |
| 0 | 18 | 19 | 0 | 0 | 22 | 23 | 0 |
| 25 | 0 | 0 | 28 | 29 | 0 | 0 | 32 |
| 33 | 0 | 0 | 36 | 37 | 0 | 0 | 40 |
| 0 | 42 | 43 | 0 | 0 | 46 | 47 | 0 |
| 0 | 50 | 51 | 0 | 0 | 54 | 55 | 0 |
| 57 | 0 | 0 | 60 | 61 | 0 | 0 | 64 |

**Note:**

i) For the main diagonal of a 4 x 4 block, i=j, where i and j represent the row and column of the matrix. For the diagonal in the opposite direction, i+j = 5.

ii) The cell value of each cell is obtained by the formula $c_{ij} = (i-1) * n + j$, where, i,j represent the row and column of the matrix and n is the order of the matrix.

| Col \ Row | 1 | 2 | 3 | 4 | 5 | 6 | 7 | 8 |
|---|---|---|---|---|---|---|---|---|
| 1 | 1 | 2 | 3 | 4 | 5 | 6 | 7 | 8 |
| 2 | 9 | 10 | 11 | 12 | 13 | 14 | 15 | 16 |
| 3 | 17 | 18 | 19 | 20 | 21 | 22 | 23 | 24 |
| 4 | 25 | 26 | 27 | 28 | 29 | 30 | 31 | 32 |
| 5 | 33 | 34 | 35 | 36 | 37 | 38 | 39 | 40 |
| 6 | 41 | 42 | 43 | 44 | 45 | 46 | 47 | 48 |
| 7 | 49 | 50 | 51 | 52 | 53 | 54 | 55 | 56 |
| 8 | 57 | 58 | 59 | 60 | 61 | 62 | 63 | 64 |

**Step 5**
Trace through each row and fill all the cells of the entire matrix with zero entries with the value of n x n – $c_{ij}$.

| 1 | 63 | 62 | 4 | 5 | 59 | 58 | 8 |
|---|---|---|---|---|---|---|---|
| 56 | 10 | 11 | 53 | 52 | 14 | 15 | 49 |
| 48 | 18 | 19 | 45 | 44 | 22 | 23 | 41 |
| 25 | 39 | 38 | 28 | 29 | 35 | 34 | 32 |
| 33 | 31 | 30 | 36 | 37 | 27 | 26 | 40 |
| 24 | 42 | 43 | 21 | 20 | 46 | 47 | 17 |
| 16 | 50 | 51 | 13 | 12 | 54 | 55 | 9 |
| 57 | 7 | 6 | 60 | 61 | 3 | 2 | 64 |

### 4.2 Matlab Code for Magic Square of order n, where n is doubly even

```matlab
% order of magic square – n=64

a=zeros(64,64);
a=double(a);
n=64; % order of magic square
for x=1:4:n
    for y = 1:4:n
        q=0;

        for i = x:x+3
            q=q+1;
            q1=0;
            for j = y: y+3
                q1=q1+1;
                if ((i==j)  || (i+j)==5)||(q+q1 == 5) || (q==q1)
                    a(i,j) = n*(i-1)+ j;
                else
                    a(i,j) = n*n-((i-1)*n+j)+1;
                end
            end
        end
    end
end
```

## 5 CONCLUSION

In this paper, a simple and easy method of constructing a doubly even magic square is proposed and implemented. The algorithm and the code in Matlab are also given so that the applications of magic squares can be explored in detail.